\documentclass[preprint,12pt,authoryear,review]{elsarticle}

\usepackage{amsmath}
\usepackage{amsfonts}
\usepackage{amssymb}
\usepackage{graphicx}
\graphicspath{}
\usepackage{natbib}
\usepackage{xcolor}

\usepackage{amssymb}

\journal{Computers \& Chemical Engineering}

\begin{document}

\begin{frontmatter}



\title{From partial data to out-of-sample \\ parameter  and observation estimation \\ with Diffusion Maps and Geometric Harmonics}


\author[inst1,inst2]{Eleni D. Koronaki}
\author[inst3]{Nikolaos Evangelou}
\author[inst3]{Yorgos M. Psarellis}
\author[inst2]{\\ Andreas G. Boudouvis}
\author[inst3]{Ioannis G. Kevrekidis\corref{cor1}}

\affiliation[inst1]{Interdisciplinary Center for Security, Reliability and Trust,
  University of Luxembourg,
  29 John F. Kennedy Avenue 1855, Luxembourg}

\affiliation[inst2]{organization={School of Chemical Engineering, National Technical University of Athens},
           addressline={9 Heroon Polytechniou str.}, 
            city={Zographos Campus},
            postcode={15780}, 
            state={Attiki},
            country={Greece}}
            
\affiliation[inst3]{
  Department of Chemical and Biomolecular Engineering,
  Whiting School of Engineering, Johns Hopkins University,
  3400 North Charles Street, Baltimore, MD 21218, USA
}

\cortext[cor1]{yannisk@jhu.edu}

\begin{abstract}
A data-driven framework is presented, that enables the prediction of quantities, either observations or  parameters, given sufficient partial data. The framework is illustrated via a computational model of the deposition of Cu in a Chemical Vapor Deposition (CVD) reactor, where the reactor pressure, the deposition temperature and feed mass flow rate are important process parameters that determine the outcome of the process. The sampled observations are high-dimensional vectors containing the outputs of a detailed CFD steady-state model of the process, i.e. the values of velocity, pressure, temperature, and species mass fractions at each point in the discretization. A machine learning workflow is presented, able to predict out-of-sample (a) observations (e.g.  mass fraction in the reactor) given process parameters (e.g. inlet temperature); (b) process parameters given observation data; and (c)  partial observations (e.g. temperature in the reactor) given other partial observations (e.g. mass fraction in the reactor).
The proposed workflow relies on the manifold learning schemes Diffusion Maps and the associated Geometric Harmonics. Diffusion Maps is used for discovering a reduced representation of the available data, and Geometric Harmonics for extending functions defined on the manifold. In our work a special use case of Geometric Harmonics is formulated and implemented, which we call Double Diffusion Maps, to map from the reduced representation back to (partial) observations and process parameters. A comparison of our manifold learning scheme to the traditional Gappy-POD approach is provided: ours can be thought of as a "Gappy DMAP" approach. The presented methodology is easily transferable to application domains beyond reactor engineering.
\end{abstract}




\end{frontmatter}




\section{Introduction}
Since nonlinear manifold learning methods were introduced \citep{balasubramanian2002isomap,roweis2000nonlinear,r19,r20,r21}, a new route was carved for the parsimonious description of data derived from models of nonlinear applications. 

The main premise of reduced order modeling methodologies is that state observables often live in low-dimensional manifolds, despite their apparent high dimensionality. Nonlinear manifold learning methodologies identify a parametrization of the manifold that describes the data \citep{xing2016manifold,r25,r14,r23,r24} and, when coupled with appropriate mapping between the reduced description and the high-dimensional ambient space, they can be used for interpolation and regression \citep{r22,evangelou2022double,evangelou2022parameter,r40}.

 Here, we demonstrate how the mapping between the ambient and the reduced space, determined with Diffusion Maps (DMAPs), can be used not only to enable efficient prediction of outputs (in our case, observations) given new inputs (in our case, process parameters), or the inputs that correspond to a new output, but also in the spirit of a static nonlinear observer \citep{kazantzis1998nonlinear,luenberger1964observing}: for the prediction of all or only part of the variables or parameters, given partial information. To achieve that, DMAPs is implemented in conjunction with a special use case of Geometric Harmonics interpolation \citep{coifman2006geometric,r40,evangelou2022double}. The latter is implemented not only as a means of mapping between the reduced and ambient space, but also as a regression tool between the input and output space. 
 
 The goal is to reconstruct variables that are inaccessible, due to technical considerations pertinent to process stability and product quality. This is particularly important in the context of dynamical systems with process control as the ultimate goal \citep{r26,r27,r28,r29,r31,r32,r33}. In this work, our methodologies are implemented for Chemical Vapour Deposition (CVD), a process where in-situ sensors are scarce and therefore critical process variables that influence product quality are inferred by ex-situ measurements. \citep{r34,GLEASON2020113,r36,r37,koronaki2014non,koronaki2016efficient,psarellis2018investigation,papavasileiou2022efficient}.
 
The proposed methodology is reminiscent of the so-called Gappy POD method, proposed by Everson and Sirovich \citep{everson1995karhunen}, as an extension of the Proper Orthogonal Decomposition (POD) method, that accounts for partly known (hence "gappy") data. According to Gappy POD, it is possible to accurately reconstruct a vector that is only partially known, provided that it is spanned by a previously defined basis of POD vectors~\citep{willcox2006unsteady,XING2022110549,JO2019419}. In the case where the data lives in a \textit{curved manifold}, the size of the POD basis required for accurate reconstruction of the data is expected to be high, since several hyperplanes are necessary to describe it; this will be discussed briefly in a subsequent section. This drawback is addressed with DMAPS, which typically require less coordinates than its linear counterpart, to accurately capture the data variance. 

 The remainder of the paper is organized as follows: the main concepts pertaining to Diffusion Maps and Geometric Harmonics are presented, followed by Double Diffusion Maps, which is the particular implementation of the latter necessary for interpolation. The Gappy POD method is then, summarized for completeness. The CVD reactor used in this work as a case study is then briefly described, as well as some details about the CFD model which generates the data. The results of the proposed workflow are then presented and compared to the Gappy POD algorithm, followed by our conclusions.

\section{Diffusion Maps}
\label{sec:Diffusion_Maps}
Diffusion maps~\citep{r19,r20,r21} is a framework based upon diffusion processes for finding meaningful geometric descriptions of data sets, even when the underlying geometry of the data is complex, nonlinear and corrupted by (relatively low level) noise. The method is based on the construction of a Markov transition probability matrix, corresponding to a random walk on a graph, whose vertices are the data points, with transition probabilities being the local similarities between data points. The first few eigenvectors of the sparse Markov matrix are used as data-driven coordinates that provide a reparametrization of the data. 

To construct a low-dimensional embedding for a data set $\mathbf{X}$ of $N$ individual points (represented as $d$-dimensional real vectors \(x_1,..., x_{N}\)), $\mathbf{X} \in \mathbb{R}^{N \times d}$ ,a similarity measure between each pair of vectors \(x_i, x_j\) is computed. The standard Euclidean distance is typically used for thid purpose. By using this similarity measure,
an affinity matrix is constructed. A popular choice is the Gaussian kernel
\begin{equation*}
w(i,j)=exp\left[-\left(\frac{||x_i-x_j||}{\epsilon}\right)^2\ \right]
\end{equation*} where $\epsilon$ is a hyperparameter that quantifies the kernel's bandwidth. 
To recover a parametrization insensitive to the sampling density, the normalization
\[\widetilde{\textbf{W}} = \textbf{P}^{-1} \textbf{W} \textbf{P}^{-1}\]
is performed, where $P_{ii} = \sum_{j=1}^N W_{ij}$ and $W_{ij}$ the elements of the matrix $\textbf{W}$.
A second normalization applied on $\widetilde{\textbf{W}}$, %
\[\textbf{K}=\textbf{D}^{-1} \widetilde{\textbf{W}}\] 
gives a \(N\times N\) Markov matrix \textbf{K}; here \(\textbf{D}\) is a diagonal matrix, collecting the row sums of matrix \(\widetilde{\textbf{W}}\).
The stochastic matrix \(\textbf{K}\) has a set of real eigenvalues \hbox{\(1=\lambda_1\geq\ ... \geq\lambda_N\)} with corresponding eigenvectors $\phi_i$. 

To check if dimensionality reduction can be achieved, the number of retained eigenvectors has to be appropriately truncated. In practice, it is useful to consider that not all obtained eigenvectors parametrize independent directions, but rather most of them can be considered as \textit{spanning the same directions with different frequencies}. Eigenvectors that parametrize the same directions in this context are called \textit{harmonics} and the ones that parametrize independent directions \textit{non-harmonics}. 
Therefore, a minimal representation of the DMAP space is made possible by carefully selecting the non-harmonic coordinates, which do not necessarily correspond to the most dominant eigenmodes of the Markov matrix. This is a stark difference between DMaps and its linear counterpart, POD (also known as Principal Components analysis), where the dominant modes are retained for the truncated representation of the data. If the number of the \textit{non-harmonic} eigenvectors is less than the number of the ambient space dimensions then model (variable) reduction is achieved.

An algorithm for identifying the non-harmonic eigenvectors is presented in~\citep{r25}, based on local linear regression. In a nutshell, a local linear function is used in order to fit the DMAP coordinate $\phi_k$ as a local linear function, \(f\), of the previous vectors $(\widetilde{\Phi}_{k-1}=[\phi_1,\phi_2,...,\phi_{k-1}])$. If $\phi_k$ can be accurately expressed as function of the other DMAP coordinates over the data, then it does not represent a new direction on the dataset and is omitted for dimensionality reduction. On the contrary if $\phi_k$ cannot be expressed as a function of the previous eigenvectors, then $\phi_k$ is a new independent eigendirection that is retained for a parsimonious representation of the data. To quantify the accuracy of the fit, the following metric is used:\[r_k=\sqrt{\frac{\sum_{i=1}^{n}(\phi_k(i)-f(\widetilde{\Phi}_{k-1}(i)))^2}{\sum_{i=1}^{n}(\phi_k(i))^2)}}\]
A small value of \(r_k\) is associated with a  $\phi_k$ that is a harmonic function of the previous eigenmodes, whereas a higher value of $r_k$ signifies that $\phi_k$ is a new independent direction on the data manifold. It has been shown in~\citep{r25} that selecting only the eigenvectors that correspond to higher values of $r_k$ leads to a parsimonious representation of the data. Eventually, the vector $x_i$ is mapped to a vector whose first component is the \textit{i}-th component of the first selected nontrivial eigenvector, whose second component is the \textit{i}-th component of the second selected nontrivial eigenvector, etc.

To map a new point, $x_{new}$, from the ambient space to DMAP space, a mathematically elegant approach known as Nyström extension is used \citep{nystrom1929praktische,fowlkes2001efficient}. The starting point of the Nyström extension is to compute the distances between the new point, $x_{new}$ in ambient space, and the $N$ data points in the original data set; the same normalizations used for DMAP need to be applied also here. The Nystr\"om extension formula reads

\[\phi_j(x_{new})=\lambda_j^{-1}\sum_{i=1}^{N}\tilde{k}(x_i,x_{new})\phi_{j}(x_{i}),\]
where $\lambda_j$ is the \textit{j}-th eigenvalue, $\phi_j(x_i)$ is the \textit{i}-th component of the j-th eigenvector and $\tilde{k}(\cdot,x_{new})$ is the kernel's value between the new point and each point in the original data set.  

\section{Geometric Harmonics}
Geometric Harmonics was introduced in \citep{r19}, inspired by the Nyström Extension, as a scheme for \textit{extending} functions defined on data $\mathbf{X}$, $f(\mathbf{X}):\mathbf{X} \to \mathbb{R}$, for $x_{new} \notin \mathbf{X}$. This extension is achieved by using a particular set of basis functions called Geometric Harmonics. These functions are computed as eigenvectors  of the symmetric \(N\times N\) \(  \textbf{W}\) matrix.  The eigendecomposition of the symmetric and positive semidefinite matrix ${\mathbf{W}}$ leads to a set of orthonormal eigenvectors $\mathbf{\psi}_1,\mathbf{\psi}_2, \dots, \mathbf{\psi}_N$ and non negative eigenvalues $\sigma_1 \geq \sigma_2 \geq \dots \geq \sigma_N \geq 0 $. 

From this set of eigenvectors to avoid numerical issues we consider a truncated subset  \hbox{\(S_\delta\) = (\(\alpha\):  \({\sigma_{\alpha}\geq\delta\sigma_1}\))} where \(\delta>0\). The extension of $f$ for a new point $x_{new}$
is accomplished by firstly projecting the function of interest in the (truncated) computed set of eigenvectors

\[f\rightarrow P_\delta f =\sum_{\alpha\in S_\delta}\langle f,\psi_\alpha \rangle \psi_\alpha,\] where $P_{\delta}$ denotes the projection of the function $f$ on the eigenvectors we retained and $\langle f,\psi_\alpha \rangle$ is the inner product between the function $f$ and the obtained $\alpha$-th eigenvector $\psi_{\alpha}$. This projection step is performed only once. 

To obtain the values of the function $f$ for $x_{new} \notin \mathbf{X}$ we extend the Geometric Harmonic Functions as

\begin{equation*}
\Psi_\alpha(x_{new})=\sigma_\alpha^{-1}\sum_{i=1}^{N}w(x_{new},x_i)\psi_{\alpha}(x_{i}) 
\end{equation*}

where $\sigma_{\alpha}$ is the $\alpha$-th eigenvalue, $\psi_\alpha(x_i)$ is the $i$-th component of the $\alpha$-th eigenvector and $w(x_{new}, x_i)$ denotes the kernel
\begin{equation*}
    w(x_{new},x_i)=exp\bigg[{-\left(\frac{||x_{new} - x_i||}{\tilde{\epsilon}}\right)^2\bigg]}
\end{equation*}
The function $f$ at $x_{new}$ is then estimated as a linear combination of the extended Geometric Harmonics
\begin{equation*}
    (Ef)(x_{new})= \sum_{\alpha \in S_{\delta}}\langle f,\psi_{\alpha} \rangle\Psi_\alpha(x_{new})
\end{equation*}
where $Ef$ denotes the estimated values of $f$ at $x_{new}$. 

\section{Double Diffusion Maps and Latent Harmonics}

A \textit{slight} twist of the  Geometric Harmonics is presented in this section. As discussed above, Geometric Harmonics constructs an input-output mapping between the ambient coordinates \textbf{X} and a function of interest $f$ defined on \textbf{X}. However, it is possible, if the data are lower dimensional to construct a map in terms of only the non-harmonic eigenvectors, \textbf{$\Phi$}. This can be achieved by operating directly on the non-harmonic DMAPs coordinates. Similar to the \textit{traditional} Geometric Harmonics, firstly an affinity matrix is constructed

\begin{equation*}
\overline{w}(i,j)=exp\left[-\left(\frac{||\phi_i-\phi_j||}{\epsilon^\star}\right)^2\ \right].
\end{equation*}
in this case the affinity matrix is constructed in term of the DMAPs coordinates. To distinguish the notation between Geometric Harmonics and Double Diffusion Maps we will use overlined symbols and $\epsilon^*$.
As in the \textit{traditional} Geometric Harmonics, the function $f$ is projected to a truncated set of the obtained eigenvectors
\[f\rightarrow P_\delta f =\sum_{\beta\in \overline{S}_\delta}\langle f,\overline{\psi}_\beta \rangle \overline{\psi}_\beta.\] 
The extension of $f$ for $\mathbf{\phi}_{new}$ is achieved by firstly extending the values of the Geometric Harmonic functions $\overline{\Psi}_{\beta}$ for $\mathbf{\phi}_{new},$

\begin{equation*}
\overline{\Psi}_\beta(\mathbf{\phi}_{new})=\overline{\sigma}_\beta^{-1}\sum_{i=1}^{N}\overline{w}(\phi_{new},\phi_i)\overline{\psi}_{\beta}(\mathbf{\phi}_{i}), 
\end{equation*}
and then estimating the value of $f$ at $\mathbf{\phi}_{new}$
\begin{equation*}
(Ef)(\mathbf{\phi}_{new})= \sum_{\beta \in \overline{S}_{\delta}}\langle f,\overline{\psi}_{\beta} \rangle\overline{\Psi}_\beta(\mathbf{\phi}_{new})
\end{equation*}.

\section{Gappy POD}

In this section the Gappy POD method is summarized, in order to better highlight the differences with the proposed approach. Lets consider the matrix, $\overline{\mathbf{X}} = \mathbf{X}^\text{T}$, the transpose of the data set $\mathbf{X}$.
A POD basis, $\mathbf{U} \in \mathbb{R}^{d \times N}$, of $\overline{\mathbf{X}}$ is computed. We approximate $\overline{\mathbf{X}}$ using a truncated number, $N^p$ of basis vectors, where $N^p \leq N$ such that:

\begin{equation*}
    \frac{\lVert  \overline{\mathbf{X}}-\widetilde{\bf{X}}\rVert} {\rVert{\overline{\mathbf{X}}}\rVert} 100 \leq \epsilon^{p} 
\end{equation*}

\noindent where $\epsilon^{p}$ is a prescribed tolerance for the truncation; in our case $\epsilon^p = 5\%$ was selected.





Let us consider a vector $x_{new}$ not in the original data set, spanned by the same basis $\bf{U}$, for which only $\mu$ values of this vector are known ($\mu$ partial observations are known) denoted as $x_{new}^{partial}$  

\begin{equation*}
    x^{partial}_{new}= m \odot x_{new}
\end{equation*}

\noindent where, $\odot$ denotes the pointwise multiplication, between the vector $x_{new}$ and $m$ the masking vector that contains ones in the positions of the $\mu$ known vector components and zeros in the rest.

The goal is to find the missing values (observations) of $x^{partial}_{new}$. This is  achieved by the following step

\begin{equation*}
    x^{rec}_{new} = \mathbf{U^{p}}c
\end{equation*}

\noindent where  $x^{rec}_{new}$ is the recovered vector, $\mathbf{U}^p$ is the truncated POD basis and $c$ are the POD coefficients estimated by solving the following linear system of equations

 \begin{equation*}
    \mathbf{A} \cdot \text{c}= (m \odot \mathbf{U}) \cdot x^{partial}_{new}
\end{equation*},

\noindent where $\mathbf{A}$ is given by
\begin{equation*}
    \textbf{A} = (m \odot \mathbf{U}) ^{\text{T}}  (m \odot \mathbf{U}) 
\end{equation*}







\subsection{The drawback of hyperplanes in parsimoniously capturing nonlinearity}
\begin{figure}[htbp]
    \centering
    \includegraphics[width=1\textwidth]{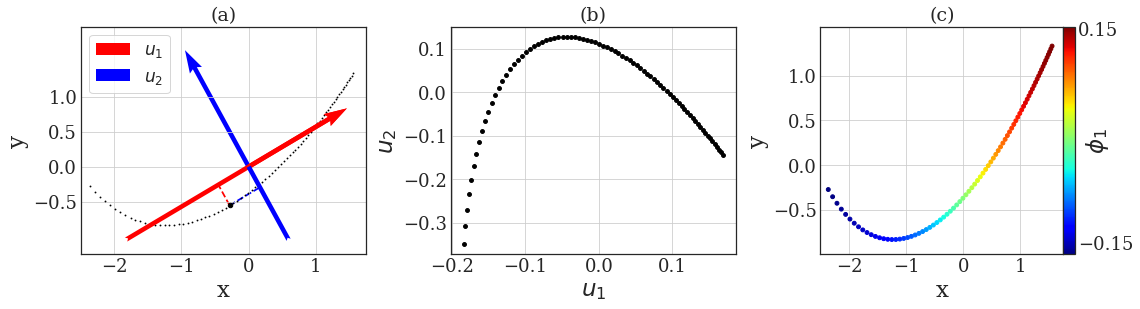}
    \caption{(a) A data set sampled from a singularly perturbed dynamical system is shown (black dots). The span of the first POD basis vector is shown with a red vector ($\textbf{u}_1$) and the span of the second POD basis vector is shown as a blue vector ($\textbf{u}_2$) (b) The components of the first POD basis vector $\textbf{u}_1$ plotted against the components of the second POD basis vector $\textbf{u}_2$ indicating that $\text{u}_2$ is a function of $\textbf{u}_1$ (c) the first non-trivial DMAPs, $\phi_1]$, eigenvector is plotted as function on the data set indicating that is able to fully parameterize it.}
    \label{fig:figure0}
\end{figure}

 In the case where the data lives in a curved manifold, the size of the POD basis required for accurate reconstruction of the data is expected to be higher than the intrinsic manifold dimension since several hyperplanes are necessary to span it. In Figure \ref{fig:figure0}, an illustrative caricature of data sampled from the singularly perturbed system $\dot{x} =  2 - x-y$, $\dot{y}= \frac{1}{\varepsilon}(x-y)$ is aiming to convey this shortcoming of POD in contrast to its nonlinear counterpart, DMAPs. In Figure \ref{fig:figure0} (a) the direction of the first POD basis vector $\textbf{u}_1$ is shown as red vector. The direction of $\textbf{u}_1$ parametrizes the data but does not fully span them because of their nonlinearity. To accuratelly reproduce this non-linear curve both POD basis vectors are needed. Please notice that, as can be seen from Figure \ref{fig:figure0} (b), the coefficient of the second POD basis vector ($\textbf{u}_2$) is a function of the first POD mode coefficient ($\textbf{u}_1$). On the contrary, Figure \ref{fig:figure0} (c) illustrates that DMAPs applied on this data set need only a single coordinate, $\phi_1$ to fully parametrize the data.

\section{Case study}
\begin{figure}[htbp]
    \centering
    \includegraphics[width=1.\textwidth]{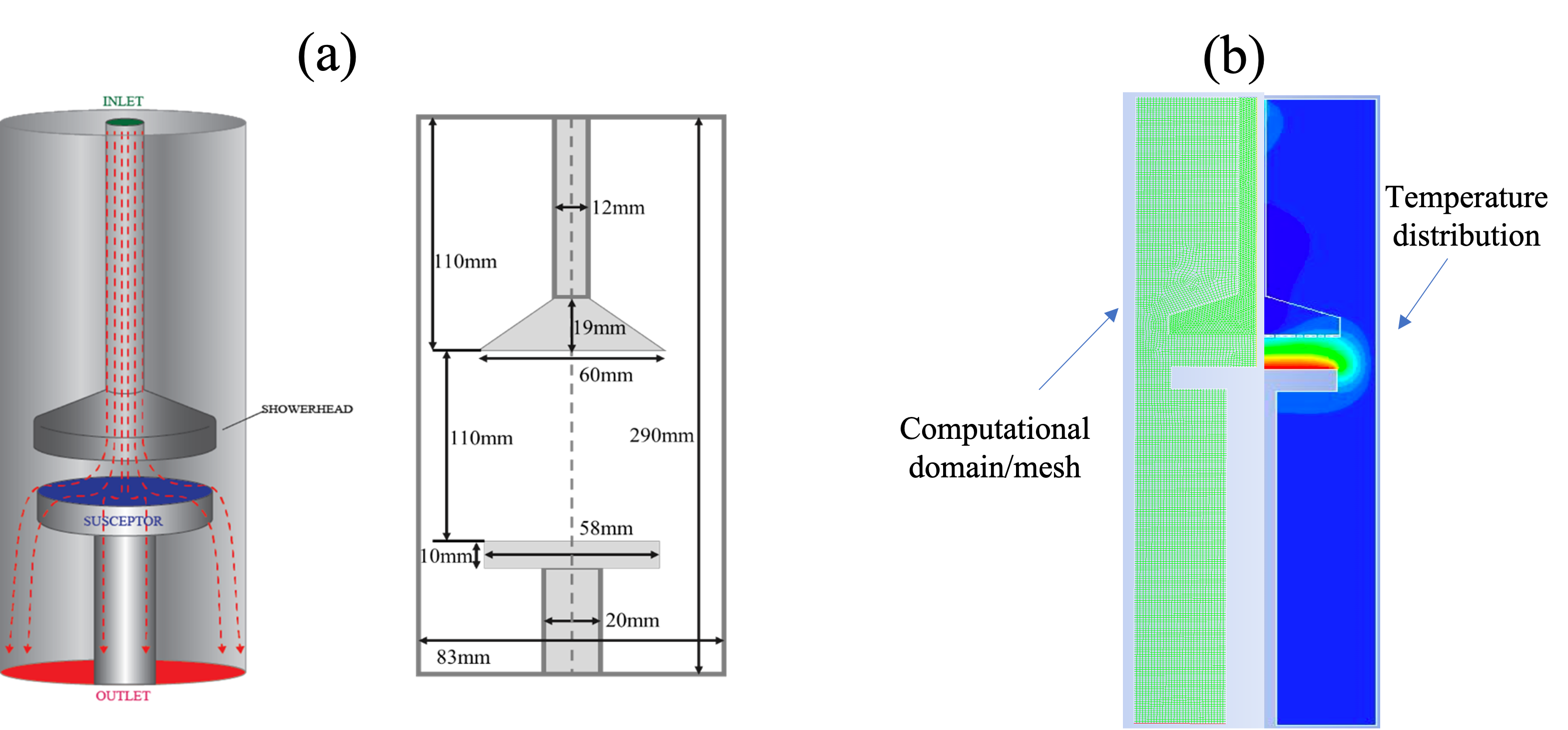}
    \caption{(a) Schematic illustration of the experimental MOCVD reactor; (b) 2D computational domain with axial symmetry and representative temperature distribution.}
    \label{fig:figure1}
\end{figure}

The case study here, is the vertical cylindrical Metal Organic Chemical Vapor Phase Deposition (MOCVD) reactor used for the deposition of Cu from copper amidinate, described in  \cite{spencer2021investigation}, shown in Figure \ref{fig:figure1}. The mixture of gas reactants enters the chamber from the top, then gets evenly distributed by passing through a showerhead and eventually leads to the deposition of Cu on a heated substrate. The quality of the produced film is affected by various process parameters; among them significant effect have the deposition temperature, $T$, the mass flow rate of incoming gas, $M$ and the chamber pressure, $P$.

\begin{figure}[htbp]
    \centering
    \includegraphics[width=1.\textwidth]{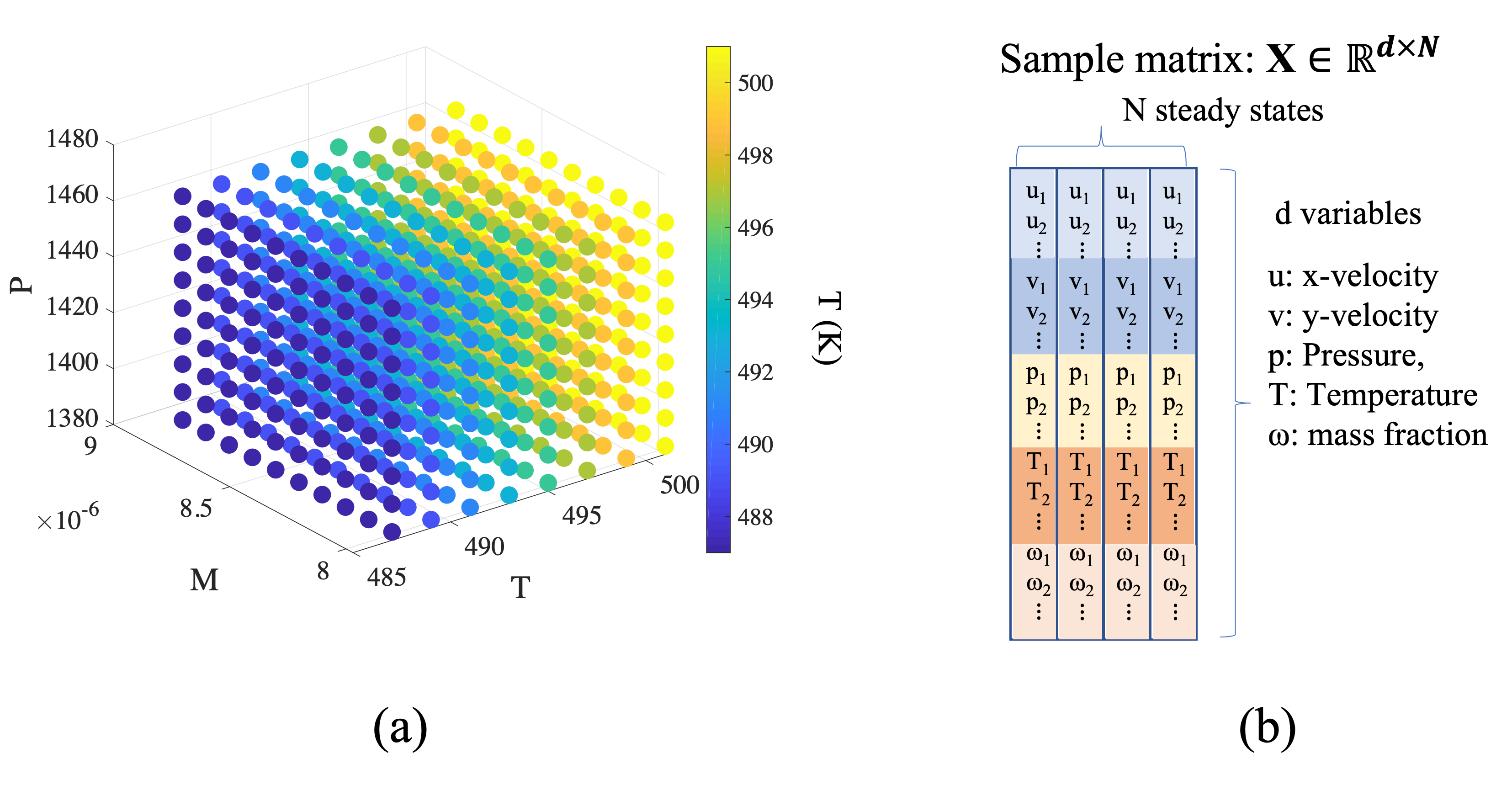}
    \caption{(a) Input parameters for the collection of data; (b) Sample set \hbox{\textbf{X} $ \in R^{(d \times N)}$}, where $ d $ is the total number of degrees of freedom of the CFD model and $N$ is the number of collected steady states.}
    \label{fig:figure2}
\end{figure}

In this work, the conservation equations for mass, momentum and energy and species are discretized with the finite volume method with 11,500 finite volumes and solved in ANSYS/Fluent, in a two-dimensional computational geometry with axial symmetry (cf. Figure \ref{fig:figure1}). The interested reader is referred to \citep{spencer2021investigation} for details on the set up of the CFD model. Here some relevant details are included for completeness.  In this implementation, the temperature of the walls and of the incoming gas mixture is constant at $T_w$ = 370 K and $T_g$=370 K respectively. The composition of the mixture of incoming gas, in terms of mass fractions is $Ar/N_2/H_2/$ Cu amidinate = 73.9\%/25.5\%/0.4\%/0.1\%. Steady states are computed for various inputs, i.e. values of three critical, for the process, parameters: the substrate temperature, T, the chamber pressure, P and the mass inflow rate of the mixture of gas reactants, $M$. Specifically, the input (or parameter) space is uniformly sampled for 487 K $<$ $T$ $<$ 501 K, 7.97 $10^{-6}$ kg/s $<$ $M$ $<$ 8.87 $10^{-6}$ kg/s and 1383 Pa $<$ $P$ $<$ 1463 Pa, as shown in Figure \ref{fig:figure2}a. This region in parameter space is interesting for the process, as it corresponds to the transition between the kinetics- and transport-limited limited regime, i.e. the process turns from being limited by the reaction rate to being defined by the diffusion rate of species toward the deposition surface.

The resulting states are collected as an ensemble of high-dimensional vectors containing the values of the two components of velocity, pressure, temperature and precursor mass fraction at each discretization point (cf. Figure \ref{fig:figure2}b). Eventually, the sample matrix $X \in R^{ N \times d}$ is assembled, where $d=58,100$ degrees of freedom (dimensions) and $N$=720 samples, i.e. vectors containing steady states.


\begin{figure}[htbp]
    \centering
    \includegraphics[width=1.\textwidth]{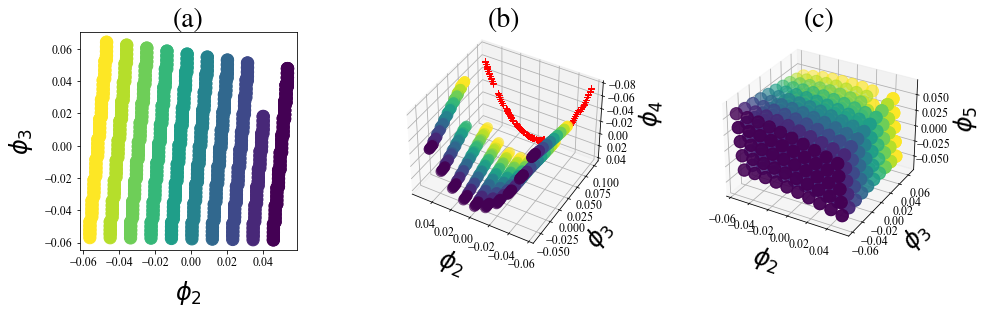}
    \caption{Diffusion Maps: (a) $\phi_3$ vs $\phi_2$; (b) $\phi_4$ vs $\phi_3$ and $\phi_2$; (c) $\phi_5$ vs $\phi_2$ and $\phi_3$; the three-dimensional ``spread" of $\phi_5$ with respect to  $\phi_3$ and $\phi_2$ suggests that these are independent directions on the low dimensional space; the distribution of $\phi_4$ with respect to $\phi_2$ and $\phi_3$ lies on a surface which indicates that $\phi_4$ is a harmonic function of $\phi_2$ and $\phi_3$}
    \label{fig:figure3}
\end{figure}

\begin{figure}[h!]
    \centering
    \includegraphics[width=.8\textwidth]{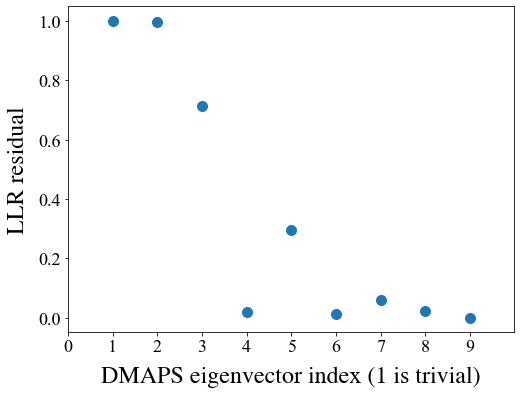}
    \caption{Residual of the local linear regression algorithm, \(r_k\); The first, second and fourth nontrivial eigenvectors \((\phi_2, \phi_3, \phi_5)\) have the highest \(r_k\) values, indicating that they each represent independent directions on the data manifold}
    \label{fig:figure4}
\end{figure}

\section{Results}
\subsection{Interpolation between ambient and intrinsic space}
The DMAPs algorithm is implemented to identify a low dimensional parametrization of the data manifold. In order to establish which are the independent coordinates, it is useful to examine the variation of each eigenvector versus the first non-trivial eigenvector of the Markov matrix, as shown in Fig. \ref{fig:figure3}: the two-dimensional variation of $\phi_3$ vs $\phi_2$ (cf. Figure \ref{fig:figure3}a) signifies that they are two independent directions on the data. Having established that the intrinsic space is at least two-dimensional, the subsequent DMAP eigenvectors are plotted against the first two. The fact the the 3D plot of $\phi_4$ vs $\phi_3$ and \(\phi_2\) reveals a surface (cf. Figure \ref{fig:figure3}b), suggests that \(\phi_4\) is a harmonic of the previous two. In contrast,  \(\phi_5\) is a new independent eigenvector and hence its variation versus the first two independent DMAPS reveals a 3D object (cf. Figure \ref{fig:figure3}c).

These visual observations are verified by the results of the implementation of the local linear regression algorithm (\cite{r25}), according to which a function $f(\phi_{k-1},\phi_{k-2},...\phi_1)$ is fitted to $\phi_k$. The results suggest that indeed \(\phi_2, \phi_3, \phi_5\) (cf. Figure \ref{fig:figure4}), represent a parsimonious low dimensional embedding of the available data, since the error, \(r_k\) of the local linear regression function is high for \(\phi_3\) and \(\phi_5\). On the contrary, \(r_k\) for \(\phi_4\) is small, indicating that it is a harmonic function of \(\phi_2\).

\begin{figure}[h!]
    \centering
    \includegraphics[width=1.\textwidth]{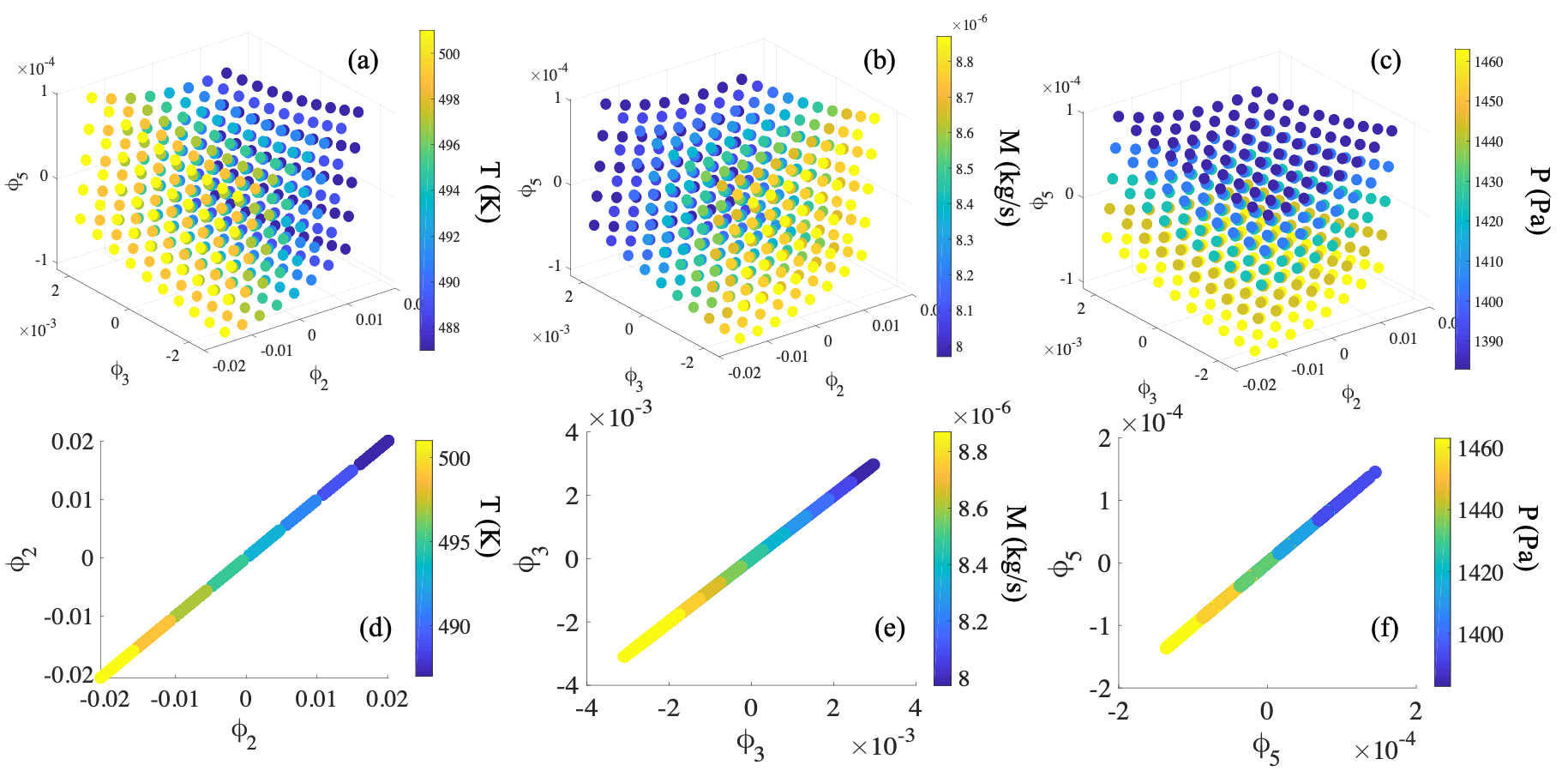}
    \caption{The identified DMAP coordinates \(\phi_2, \phi_3, \phi_5\), plotted and colored by the process inputs (a) $T$, (b) $M$, (c) $P$; the three-dimensional plots show that the variation of \(\phi_2, \phi_3, \phi_5\), follow the variation of $T$, $M$ and $P$ respectively. This is further demonstrated in (d), (e) and (f) respectively}
    \label{fig:figure5}
\end{figure}

In an attempt to provide a physical interpretation of the derived low dimensional coordinates, the eigenvectors \(\phi_2, \phi_3, \phi_5\) are plotted and colored by the values of $T$ (Figure \ref{fig:figure5}a), $M$ (Figure \ref{fig:figure5}b) and $P$ (Figure \ref{fig:figure5}c). Each one of the three directions in the reduced space corresponds to the variation of each one of the three input parameters. This is shown with more clarity in (Figure \ref{fig:figure5}d), (Figure \ref{fig:figure5}e) and (Figure \ref{fig:figure5}f), specifically that \(\phi_2\) corresponds to T, \(\phi_3\) to M and \(\phi_5\) to P. Therefore, the DMAP coordinates provide a parametrization that appears to be one-to-one with the actual physical parameters that were varied in order to produce it. 

\begin{figure}[h!]
    \centering
    \includegraphics[width=.8\textwidth]{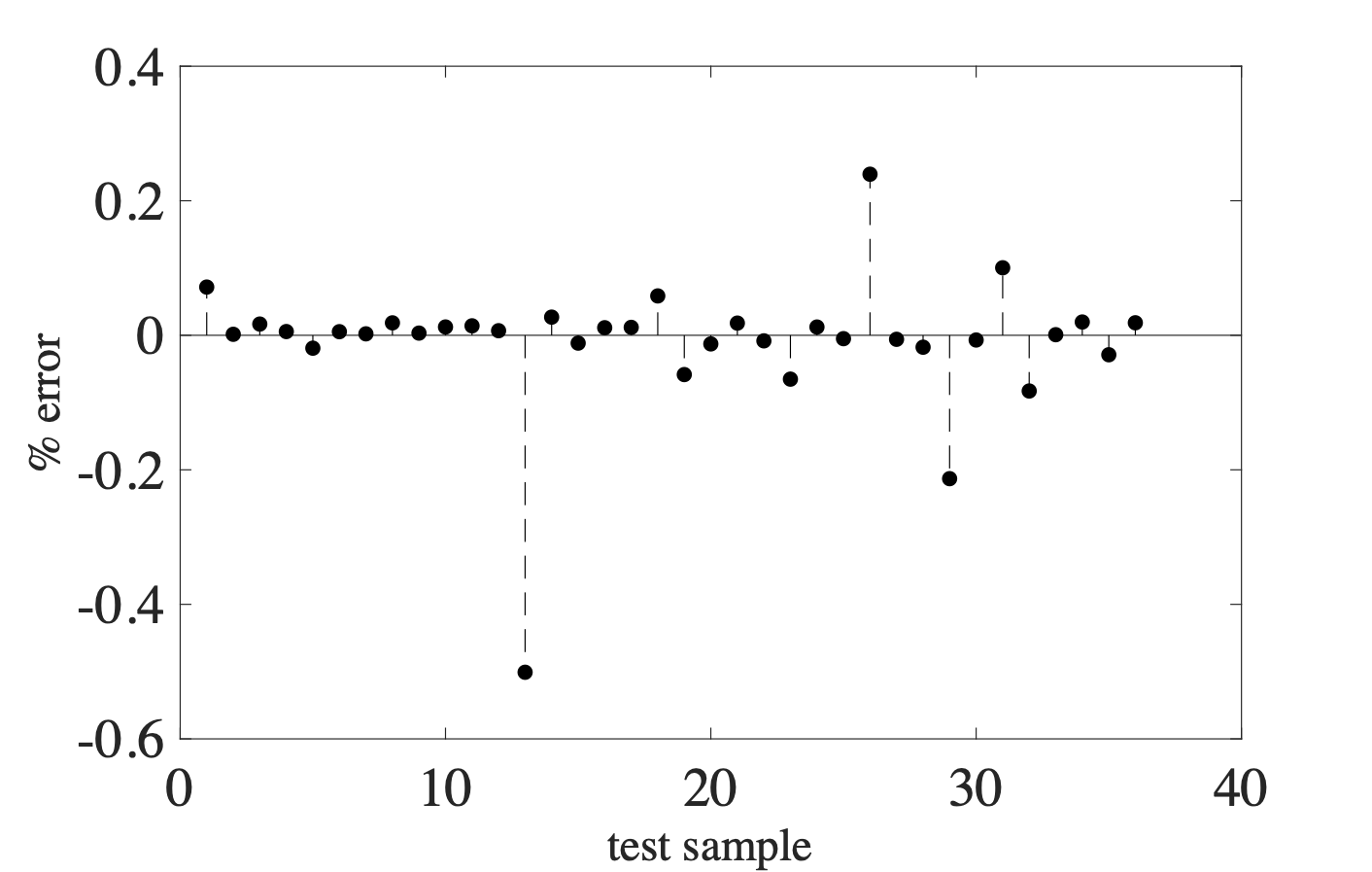}
    \caption{Prediction error of the Double DMAP implementation of Geometric Harmonics for each one of the vectors $x_i$ in the test sample}
    \label{fig:figure6}
\end{figure}

The inverse map, \(f^{-1}\), i.e. from the reduced to the ambient space, is approximated with the double DMAPS Geometric Harmonics interpolation and its accuracy is assessed against a random test sample. For this implementation, the value of the kernel parameter \(\epsilon\) is 5.105 and 38 eigenvectors are retained as interpolation functions. The \(\%\) relative error for each vector, shown in Figure 6, is computed as:\[\% error=\left(\frac{X_i^{predicted}-X_i^{actual}}{X_i^{predicted}}\right)100\]
The average error for the test samples is 0.01 \(\%\).

\subsection{Prediction of outputs for a new set of inputs}
\begin{figure}[htbp]
    \centering
    \includegraphics[width=1.\textwidth]{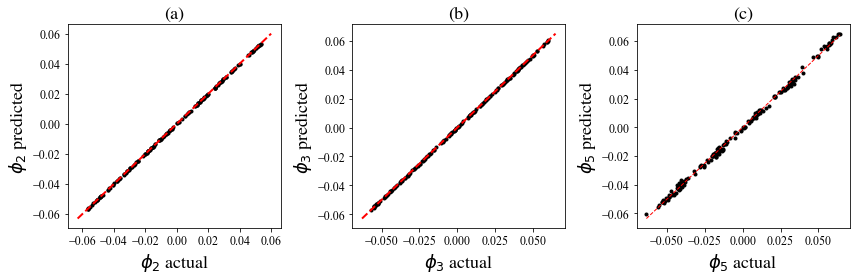}
    \caption{The predicted versus the actual DMAP coordinates. The solid lines correspond to \(y=x\).}
    \label{fig:figure7}
\end{figure}

The mapping from the ambient space to the intrinsic space and back provides a means of deriving various useful correlations, one them being the prediction of the output of a new set of input parameters. Specifically, Geometric Harmonics interpolation is implemented in order to define functions \(\phi_2 = g_1(T,P,M), g_2: (\phi_3) = g_2 (T,P,M)\) and \(g_3: (\phi_5) = g_3 (T,P,M)\). The interpolation mean squared error for \(\phi_2, \phi_3\) and \(\phi_5\) is \(7.39 10^{-8}\), \(2.23 10^{-6}\) and \(9.12 10^{-6}\) respectively. The predicted reduced coordinates versus the actual ones are shown in Figure \ref{fig:figure7}.

Having established a mapping from the input space to the reduced variables Double DMAPs interpolation with Geometric Harmonics whicjj defines the inverse map described in the previous paragraph, is implemented in order to find the corresponding state variables in the ambient space. The average mean squared error for the test sample is $(1.2 10^{-5})$.

\subsection{Prediction of the inputs that correspond to a new output}
\begin{figure}[!]
    \centering
    \includegraphics[width=1\textwidth]{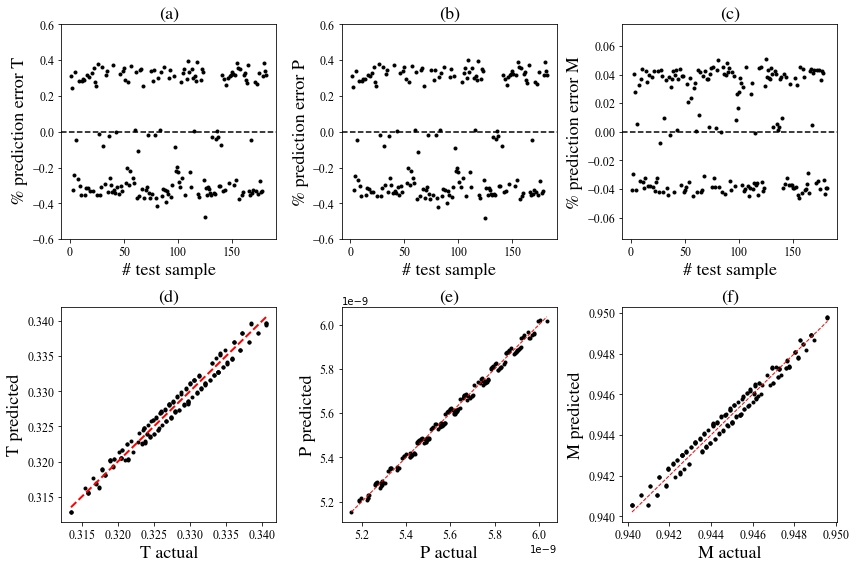}
    \caption{Prediction \% relative error of input parameters (a) $T$, (b) $P$, (c) $M$ that correspond to a new state in the ambient space. Predicted versus actual values of (d) $T$, (e) $P$, (f) $M$; the red dashed lines correspond to \(y=x\).}
    \label{fig:figure8}
\end{figure}

It is possible to find the values of inputs, ($T$, $P$, $M$) that correspond to a new output, by first using the Nystr\"om extension to obtain the corresponding reduced variables. Then an interpolation function  can be constructed from the reduced coordinates to the input parameters, with Double DMAPS: \((T,P,M)=f(\phi_2,\phi_3,\phi_5)\). The interpolation relative error is below 0.5\(\%\), as shown in the top row of Figure  \ref{fig:figure8} and the predicted input parameters versus the actual ones are shown in the bottom row of Figure \ref{fig:figure8}.

\subsection{Prediction of inputs from partial observations}
Instead of a full state vector in ambient space, it is also possible to use partial observations, such as, for example, the value of temperature at a few points, in order to predict the corresponding input values. To this end, Geometric Harmonics interpolation is implemented in order to define a function \((\phi_2,\phi_3,\phi_5)=f(X_{partial}^1,X_{partial}^2,X_{partial}^3,...,X_{partial}^m)\).  In this implementation, \(X_{partial}^i, i=1,..., m\) are values of temperature in seven random positions (cf. Figure \ref{fig:figure9}a). From the reduced space, it is now possible to map to the input space with the interpolation function from the reduced space to the input space discussed in the previous paragraph.

\begin{center}
\begin{figure}[!]
    \centering
    \includegraphics[width=1.1\textwidth]{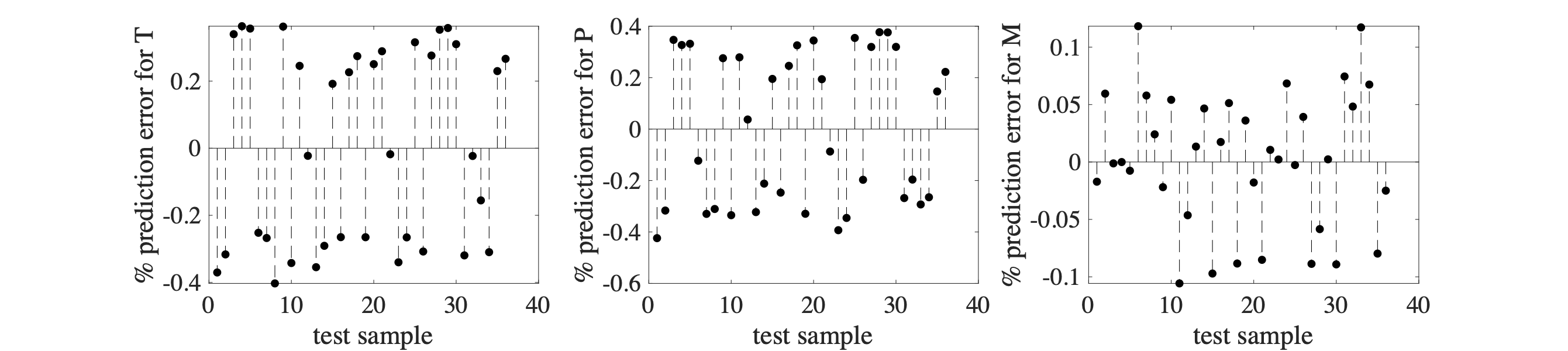}
    \caption{Prediction error of inputs $T$, $P$, $M$.}
    \label{fig:figure10}
\end{figure}
\end{center}

\begin{figure}[!]
    \centering
    \includegraphics[width=1.\textwidth]{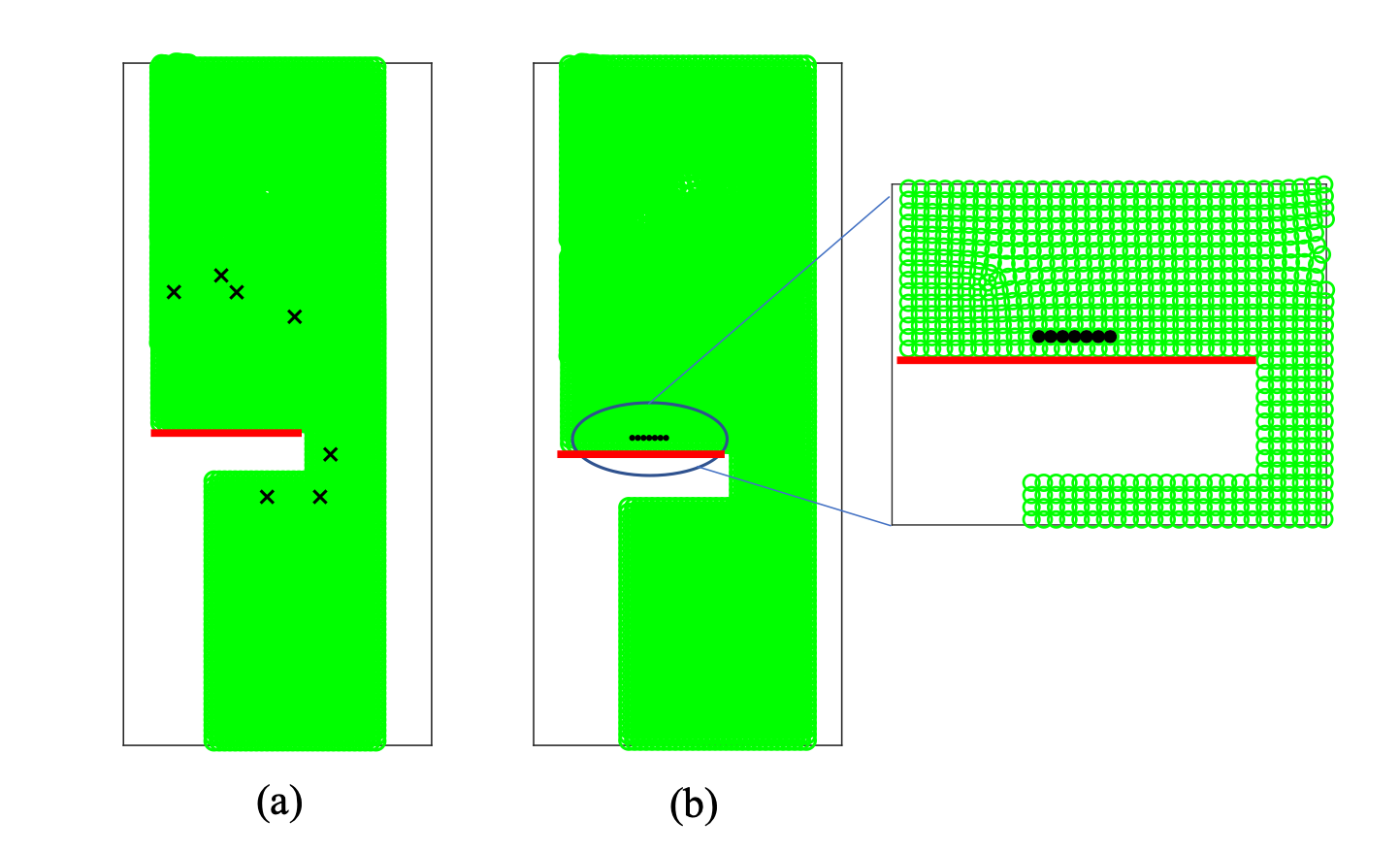}
    \caption{Partial observations: (a) positions in computational geometry where temperature values are considered; (b) positions in computational geometry where mass fractions are predicted; The red line indicates the heater susceptor surface where deposition of material takes place.}
    \label{fig:figure9}
\end{figure}

The number of partial observations required in order to define the function from the partial observations to the input space  is dictated by Whitney’s embedding theorem (\cite{whitney1936differentiable}). $2n+1$ independent observations are provably sufficient to create an embedding of the m-dimensional manifold. Here, for the three-dimensional reduced manifold, at least seven partial observations should be considered. Eventually the predicted values of the input parameters for the test sample are in excellent agreement with the actual values, as presented in Figure \ref{fig:figure10}, where the prediction relative error for the test sample is shown to be on average less than 0.1 \(\%\).

\subsection{Prediction of partial observations from  other partial observations}
The possibility to predict part of the observations, given a different part of the observations, is discussed in this paragraph. As an illustrative example, the prediction of the value of mass fractions right above the heated susceptor surface (cf. Figure \ref{fig:figure9}b), given seven temperature measurements in a different part of the geometry (cf. Figure \ref{fig:figure9}a) will be presented here. This choice is dictated by the fact that, although in this particular process the mass fraction of precursor reaching the deposition surface is crucial for determining both the quality of the product and also the film deposition rate, it is not easily measurable. On the other hand, temperature measurements are generally more accessible, and the idea is to use such measurements in order to make predictions for quantities that are harder to measure, similar to the concept of a nonlinear observer.

\begin{figure}[!]
    \centering
    \includegraphics[width=1.\textwidth]{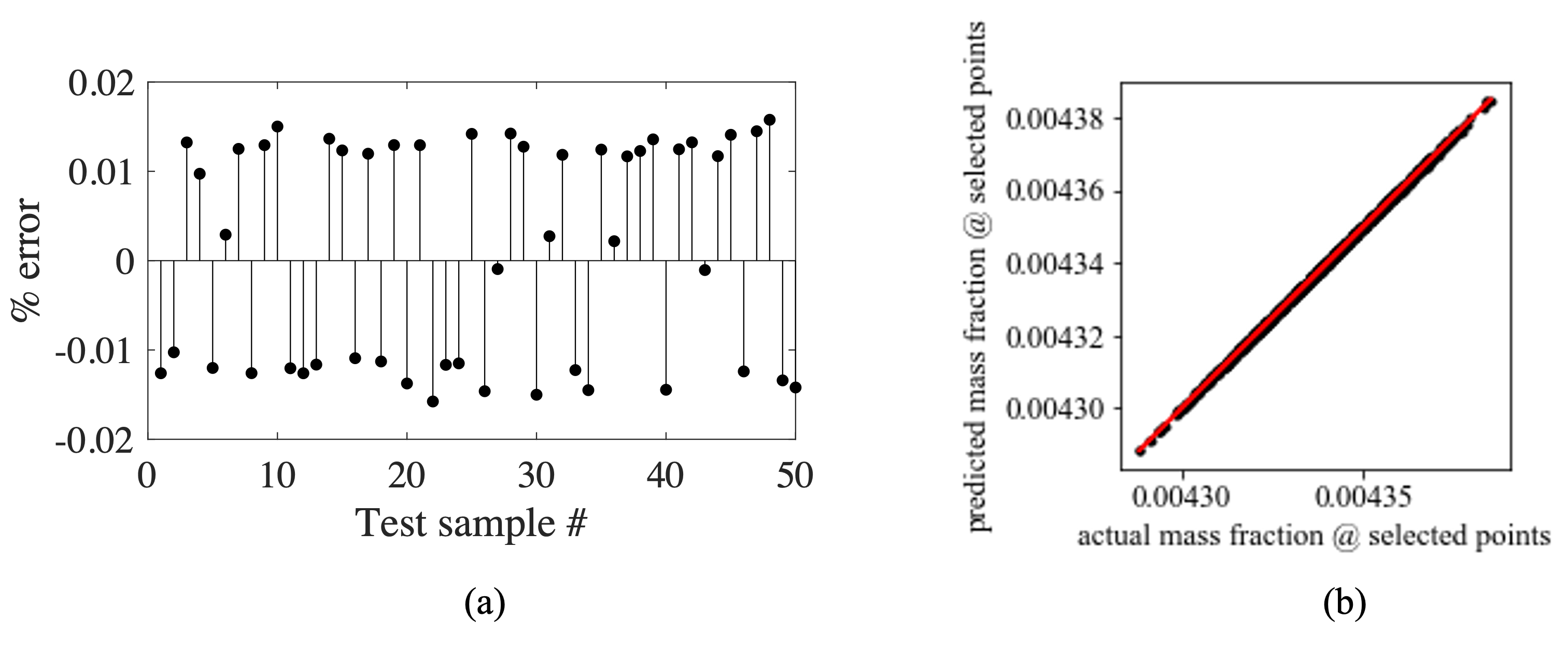}
    \caption{(a)Approximation error for the mass fraction, $\omega$ values in the sample test set (b) predicted vs actual mass fraction values at one point for all test samples}
    \label{fig:figure11}
\end{figure}

To begin with, the function $(\phi_2, \phi_3, \phi_5)=f(X_{partial}^1,X_{partial}^2,\dots,X_{partial}^m)$, discussed in the previous paragraph is used in order to map from the partial observations to the reduced space. The inverse function, \(f^{-1}\), from the reduced to ambient space, is then used in order to predict the desired values( in this case the values of the species mass fractions above the heated substrate at seven points). The average relative error is less than 0.5 \(\%\) (cf. Figure\ref{fig:figure11}a), whereas the predicted versus the actual values of the mass fraction at a single point above the substrate is shown in  Figure\ref{fig:figure11}b.

\subsection{Implementation of Gappy POD and comparison to DMAP-based predictions}

In this section the Gappy POD method is implemented and the results are compared to those delivered by the DMAPS/Geometric Harmonics method presented in the previous sections.

\begin{figure}[h!]
    \centering
    \includegraphics[width=1.\textwidth]{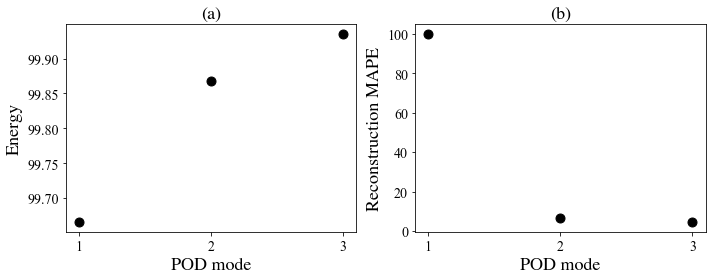}
    \caption{(a)Cumulative energy captured by POD modes (b) Mean absolute percentage error of the reconstructed data-set $\bf{X}$}
    \label{fig:figure12}
\end{figure}

The first step is to compute a basis of the given data-set using singular value decomposition. The size of the basis is determined based on the cumulative percentage of the energy of the system captured by $i$ modes, defined as
\begin{equation*}
  E_i \% = \frac{\sum_{n=1}^{i} \xi_n}{\sum_{n=1}^{m} \xi_n} *100
\end{equation*}

where $\xi_n$ stands for the $nth$ singular value of the diagonal matrix $\mathbf{\Xi}$ that results from the singular value decomposition of the transpose of the data matrix $\mathbf{\overline{X}}$. In addition to that, the reconstruction error of the data-set is computed for increasing size of the POD basis. In this implementation, and for the purposes of comparison, the selected POD basis consists of 3 vectors that capture 99.93\% of the energy of the system (cf Figure\ref{fig:figure12}a) and the approximation error is 4.7\% (cf Figure \ref{fig:figure12}b).



\begin{figure}[h!]
    \centering
    \includegraphics[width=1.\textwidth]{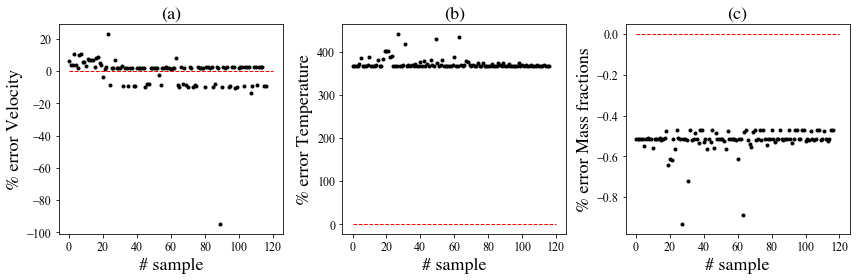}
    \caption{Prediction error (a) x-velocity  (b) temperature  (c) precursor mass fraction}
    \label{fig:figureGappy}
\end{figure}

Initially, the goal is to predict the output vector $x_{new}$, containing the distributions of velocity, pressure, temperature and mass fractions, given a new set of process parameters ($T$, $P$, $M$). Therefore, the partial data considered in this comparison correspond to the values of the three process parameters. The predicted values are compared to the projection of the test vector on the POD basis. In this case the predictions are inaccurate, especially for the temperature distributions (cf. Figure \ref{fig:figureGappy}b) and even unphysical as negative values for the mass fractions are produced (cf. Figure\ref{fig:figureGappy}c).  

\begin{figure}[h!]
    \centering
    \includegraphics[width=1.\textwidth]{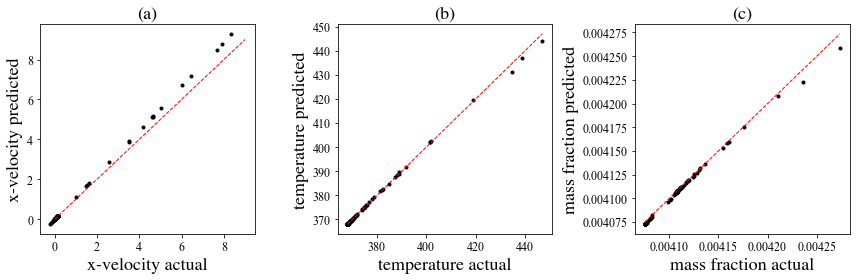}
    \caption{Predicted vs actual values of (a) x-velocity  (b) temperature  (c) precursor mass fraction}
    \label{fig:figure13}
\end{figure}

Nevertheless, the results of the Gappy POD method are heavily influenced by the choice of the ``known" values, i.e. the partial data considered. To illustrate this, instead of the three process parameters, three values of temperature are considered known, a subset of the values mentioned in Section 7.4 and shown in Figure \ref{fig:figure9}a. In this case the approximation error drops to 8\%, while the predicted versus actual velocity, temperature and mass fractions are shown in Figure \ref{fig:figure13}a, \ref{fig:figure13}b and \ref{fig:figure13}c respectively. 

This finding is directly related to the condition number of the matrix $\mathbf{A}$, defined in Section 5. Specifically, the elements of this matrix result from the inner products of the "gappy" POD vectors, i.e. the elements of the original POD vectors that correspond to the known elements of $x_{new}$. These are no longer orthogonal, and therefore the matrix $\mathbf{A}$ is fully populated. In general, the positions of the known elements, and hence the non-zero elements of  $\mathbf{A}$, must be such that orthogonality is preserved. Furthermore, the diagonal entries of $\mathbf{A}$  must not be very small, which means that the POD basis element at that point must not be small. These two requirements are reflected in the condition number of the matrix  $\mathbf{A}$: specifically, the smaller the condition number, the more they are satisfied. This analysis is presented in \cite{willcox2006unsteady}, in the context of optimal sensor placement, and in \cite{alonso2004optimal,alonso2004optimal2}, where the angle between the measurement subspace and the low dimensional space that spans the data is taken into account. The condition number of  $\mathbf{A}$ drops from $12.4e+12$ to $24.5$ when the known components correspond to the input parameters and the temperature measurements respectively. 

\begin{figure}[h!]
    \centering
    \includegraphics[width=1.\textwidth]{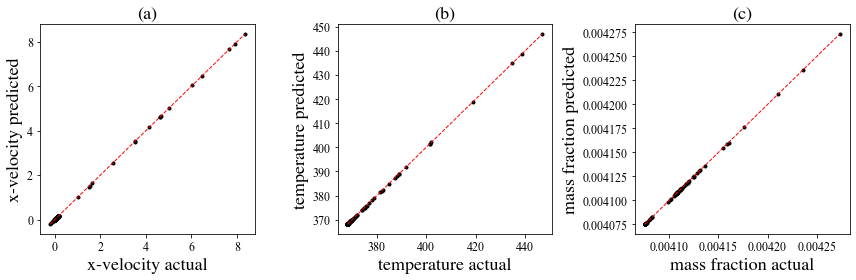}
    \caption{Predicted vs actual values of (a) x-velocity (b) temperature (c) precursor mass fractions}
    \label{fig:figure14}
\end{figure}

\begin{figure}[h!]
    \centering
    \includegraphics[width=1.\textwidth]{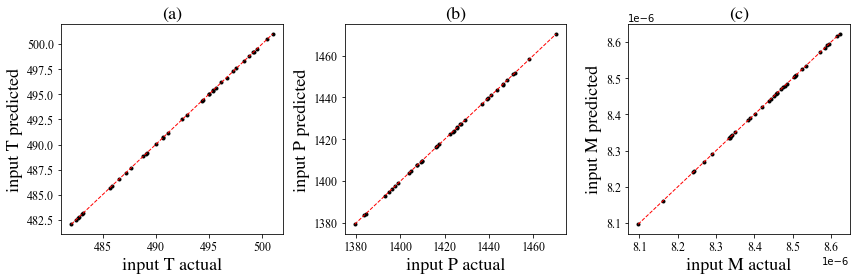}
    \caption{Predicted vs actual values of (a) temperature, $T$ (b) mass flow rate, $M$  (c) pressure, $P$}
    \label{fig:figure15}
\end{figure}

When all the temperature measurements at points shown in Figure \ref{fig:figure10}a are known, the prediction error drops further to $0.1\%$ and it is possible to reproduce accurately the distributions of velocity (cf Figure \ref{fig:figure14}a), temperature (cf Figure \ref{fig:figure14}b) and pressure (cf Figure \ref{fig:figure14}c), as well as the corresponding process parameters (cf Figure \ref{fig:figure15} a, b and for $T$, $P$ and $M$ respectively). In this case the condition number of the matrix $\mathbf{A}$ is $12.45$.

The results above point directly to the apparent disadvantage of Gappy POD, when compared to the proposed methodology, based on DMAPS: given the same number of POD vectors as DMAP coordinates, the accuracy of Gappy POD is inherently linked to which elements of the partial vector are known. No such consideration need be paid when Diffusion Maps/Geometric Harmonics are implemented, which enables the accurate prediction of entire vectors of outputs for various new combinations of input parameters, as well as from partial measurements.

\begin{figure}[h!]
    \centering
    \includegraphics[width=1.\textwidth]{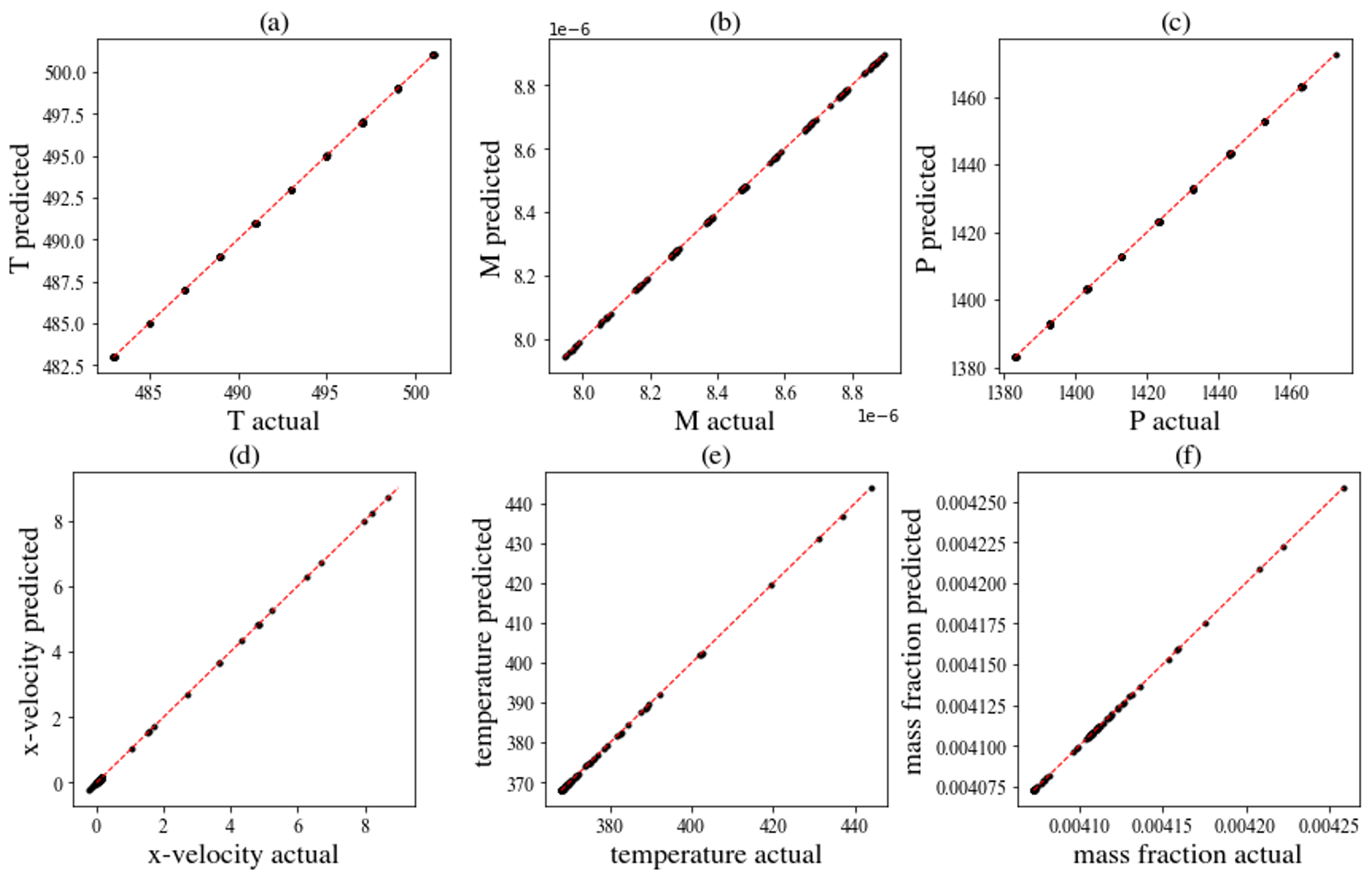}
    \caption{Gappy POD performance: 5 POD vectors and 5 known Temperature values; Predicted vs actual values of (a) temperature, $T$ (b) mass flow rate, $M$  (c) pressure, $P$, (d) x-velocity (e) temperature (f) mass fractions}
    \label{fig:figure18}
\end{figure}

Apart from the pathology related to the known values, one expected drawback of Gappy POD is directly linked to the inability of hyperplanes to accurately parametrize a curved manifold. In this implementation this is reflected in the size of the POD basis required to reconstruct the data set with an error of less than $1 \%$: here 5 POD modes are required to achieve $0.7 \%$ reconstruction error, versus the, only, 3 DMAP coordinates that are sufficient to parametrize the manifold. By selecting 5 POD modes, it is no longer possible to reconstruct the state vector, given the values of the \textit{three} input parameters, $T, P$ and $M$. In this case the values of at least 5 measurements are necessary to achieve accurate reconstruction of the state vector. This is illustrated in Figure \ref{fig:figure18}, where the predicted input parameters and distributions of velocity, temperature and mass fractions are plotted against the actual values. In this implementation, the value of temperature at five positions are considered known.

\section{Conclusions}
This work presents a data-driven workflow, based on nonlinear manifold learning, specifically Diffusion Maps, that enables the parsimonious description of high-dimensional data, but also interpolation and regression for out-of-sample predictions with remarkable accuracy. 

The case study here is a Chemical Vapor Deposition Reactor, although the proposed approach is not restricted to a particular application. Mapping between the reduced (or DMAP) and the ambient space is achieved with Geometric Harmonics, amended with a special "twist" that implements a second round of Diffusion Maps in order to define functions for accurate interpolation.

Having defined the reduced description of the data-set and the means to map back and forth between ambient and reduced DMAP space, we proceed to show the implementation of out-of-sample predictions: We first demonstrate the possibility to predict the high-dimensional output of a new set of inputs, here process parameters, namely temperature, pressure and the mass inflow rate, without additional expensive CFD simulations. The opposite, i.e. accurately predicting the inputs that correspond to a new output, is also possible.

Based on the reduced description of the data-set, provided by Diffusion Maps and the computational 
means of transitioning between the ambient and the reduced DMAP space of the data, we show how to predict not only outputs but also inputs, i.e. process parameters, when only a handful of measurements, temperature in this case are known we also demonstrate the superiority of interpolating on nonlinear manifolds (our "Gappy DMAP" approach) rather than on linear hyper-planes, as proposed by Gappy POD.

\section*{Acknowledgements}
The work of YGK, NE and YMP was partially supported by the US AFOSR and the US DOE. E.D.K. has received funding from the European Union’s Horizon 2020 research and innovation programme under the Marie Skłodowska-Curie grant agreement No 890676  - DataProMat "


\bibliographystyle{elsarticle-harv} 
\bibliography{ekor_references}





\end{document}